\documentclass[12pt,reqno]{amsart}   	
\usepackage[letterpaper, margin=1in]{geometry}
\usepackage{amssymb}
\usepackage{amsmath}
\usepackage{amsfonts}
\usepackage[utf8]{inputenc}
\usepackage{amsthm}
\usepackage{tikz-cd}

\newtheorem{thm}{Theorem}[section]
\newtheorem{lemma}[thm]{Lemma}
\newtheorem{df}[thm]{Definition}
\newtheorem{prop}[thm]{Proposition}
\newtheorem{ex}[thm]{Example}
\newtheorem{cor}[thm]{Corollary}

\theoremstyle{obs}
\newtheorem{obs}[thm]{Observation}

\theoremstyle{prob}
\newtheorem{prob}[thm]{Problem}
\theoremstyle{conj}
\newtheorem{conj}[thm]{Conjecture}

\newtheoremstyle{remark}
    {\dimexpr\topsep/2\relax} 
    {\dimexpr\topsep/2\relax} 
    {}          
    {}          
    {\bfseries} 
    {.}         
    {.5em}      
    {}          

\theoremstyle{remark}
\newtheorem{remark}[thm]{Remark}

\newcommand{\Z}{\mathbb{Z}}

\newcommand{\R}{\mathbb{R}}
\newcommand{\Q}{\mathbb{Q}}
\newcommand{\N}{\mathbb{N}}
\newcommand{\OO}{\mathcal{O}}

\title{Sums of squares on hypersurfaces}
\author{Kacper Błachut and Tomasz  Kowalczyk}
\date{}

\begin{document}

\maketitle

\begin{abstract}

We show that the Pythagoras number of rings of type $\R[x,y, \sqrt{f(x,y)}]$ is infinite, provided that the polynomial $f(x,y)$ satisfies some mild conditions.
\end{abstract}

\section*{Introduction}

One of the important problems of real algebraic geometry is the question of sums of squares. Let $f$ be a nonnegative function on some topological space $X$: is $f$ a sum of squares of functions? if so, how many squares are needed to represent $f$? This problem appears in various settings and instances, cf. \cite{Banecki, benoist, benoist 2, Cassels 2, CLDR, CLR, CLRR, fernando, to, to2}.  Here we will be interested in the latter problem in the case of polynomial functions on an algebraic hypersurface in $\R^3$.

\begin{df}
For a commutative ring $R$ with identity its Pythagoras number, $p(R)$ is the smallest positive integer $g$ such that any element which is a sum of squares can be written as a sum of at most $g$ squares. If such number does not exist, we put $p(R)=\infty$.
\end{df}
One of the most famous problems in this vein was the Hilbert 17th problem: does every nonnegative polynomial on $\R^n$ is a sum of squares of rational functions? It was solved in affirmative by Artin \cite{Artin} in 1927. Then came the first explicit example of a non-negative polynomial that is a sum of squares of rational functions, but not polynomials, namely $M(x,y)=x^2y^4+x^4y^2 -3x^2y^2+1$, given by Motzkin in 1967 \cite{Motzkin}. Not much was known about the possible number of summands until the work of Pfister \cite{Pfister} who showed that $p(\R(x_1,\dots,x_n))\leq 2^n$. Currently, the only known general bounds for the Pythagoras number of the field of rational functions $\R(x_1,x_2,\dots, x_n)$ are of the form
$$n+2 \leq p(\R(x_1,x_2,\dots, x_n))\leq 2^n$$
for $n\geq 2$ (see \cite{CLR}) which in particular gives $p(\R(x,y))=4$. Nothing more is known for $n>2$.

For arbitrarily field $K$, the Cassels Theorem \cite{Cassels} shows $p(K(x))=p(K[x])$. It is known that $p(\R(x))=p(\R[x])=2$ and $p(\Q(x))=5$ (see \cite{Pourchet}). Hoffman proved that for any positive integer $s$ there exists a real field $K $ such that $p(K)=s$ (see \cite{hoffman}). Currently, it is still not known whether $p(K)<\infty $ implies $p(K(x))<\infty$, it does hold for some special classes of fields, but the general answer is still elusive.

The Pythagoras number of a ring is a much more subtle problem. Affine $k$-algebras of transcendence degree 1, where $k$ is a real closed field, always have finite Pythagoras number. It is also known that $p(\Z[x])=\infty$ as well as $p(\R[x_1,\dots, x_n])=\infty$ for $n\geq 2$, but $p(\R[[x,y]])=2$ and $p(\Z[[x]])=5$. On the other hand, any formally real affine $k$-algebra of dimension at least three have infinite Pythagoras number, for a field $k$. These results were proven in \cite{CLDR}. Quite a lot is known for local rings, especially for local factorial rings of dimension 2 (cf. \cite{Scheiderer}). For some recent developments in this direction, see \cite{Miska} for the case of local Henselian rings, and \cite{Banecki, Banecki 2} for the rings of $k$-regulous functions. Let us briefly note that one can attempt to parametrize the space of all presentations of a given polynomial as a sum of squares, for this problem see \cite{Scheiderer 2, Vill, Vill 2} and references. This discussion suggests, that the case of dimension 2 is the most interesting case.

There are very few examples of computations of the Pythagoras number of a coordinate ring $\R[X]$ of an algebraic surface $X$. Besides $X=\R^2$ and $X=\{x^2+y^2-z^2=0 \}$ the other known examples is the family $\R[x,y,z]/(x^a-y^bz^c)$
where $a,b,c$ are arbitrary nonnegative integers, not all of them zero, each such ring has infinite Pythagoras number. Note that in the above family most surfaces are singular. All of the above examples come from the paper \cite{CLDR}.

In this paper we consider the problem of computing the Pythagoras number of a different family of coordinate rings of real algebraic surfaces in $\R^3$, surfaces defined by $\{z^2-f(x,y)=0 \}$ with some mild conditions on the polynomial $f(x,y)$. In particular, we are able to compute the Pythagoras number of coordinate rings of a large family of hypersurfaces, both singular and nonsingular. All of the above Pythagoras numbers turns out to be infinite.

 In the first section we define an admissible polynomial and recall facts and definitions regarding the  Pythagoras number of a ring. We then proceed to prove the main result of the paper in Section 2. The last section contains application of the Theorem \ref{main} to the family of hypersurfaces with the du Val singularities. We then propose some open problems and possible further research.

\section{Preliminaries}
Let us recall some definitions. 

\begin{df}
Let $R$ be a commutative ring with identity. We define the length of an element $a \in R$ to be the smallest positive integer $g$ such that $a$ can be written as a sum of at most $g$ squares.
\end{df}
If an element is not a sum of any number of squares, then its length is infinite. In general, only elements of finite length will be interesting.

\begin{df}
Let $f(x,y) \in \R[x,y]$ be a non-constant polynomial such that $f(0,0)=0$. Let $d = \deg_y f(x,y)$ and $\alpha x^b y^d$ be a monomial with the largest possible $b$ and nonzero coefficient $\alpha$ among all monomials in $f(x,y)$.
We say that $f(x,y)$ is strictly admissible if:
\begin{itemize}
\item[a)]$\alpha>0$ and $b,d$ are even, or
\item[b)] $b$ or $d$ is odd.
\end{itemize}
\end{df}

\begin{df}
We say that a polynomial $f(x,y)$ is admissible, if there exists an invertible matrix $M \in \text{GL}_2(\R)$ such that $f\circ M$ is a strictly admissible polynomial.

\end{df}

\begin{ex}
Among polynomials $f_1=x^3y$, $f_2=x^2-y^2$, $f_3=-y^2-x^7$, $f_4=-2x^2-3x^4y^2$, only $f_1$ is strictly admissible while $f_2$ and $f_3$ are admissible. 

\end{ex}

\begin{ex}
Consider the polynomial $f(x,y)=-x^2y^4-x^4y^2+3x^3y^3$, clearly, it is not strictly admissible. However, after applying a linear change of variables given by matrix
$\begin{pmatrix}
1 & -1 \\
1 & 1 
\end{pmatrix}$ we obtain $f\circ M(x,y)=f(x-y,x+y)=(x-y)^2(x+y)^2(5y^2-x^2)$ which is strictly admissible.
\end{ex}

From the very definition of admissibility, we see that an admissible polynomial is strictly positive on a sufficiently large subset of $\R^2$. One can ask, whether a polynomial which is positive on an unbounded set is admissible. This is not the case

\begin{ex}\label{example}
Let  $g(x,y)=x^2(1-x^2)$. Such a polynomial is strictly positive on an infinite strip of width 2. After an invertible linear change of coordinates we obtain $g(ax+by,cx+dy)=(ax+by)^2(1-(ax+by)^2)$. For such polynomial, the coefficient at $y^4$ is either zero or negative, hence this polynomial cannot be strictly admissible.

\end{ex}

Let us now recall some facts concerning length and Pythagoras number. Let $R_1,R_2$ be commutative rings with identity. We recall a property of the Pythagoras number that follow immediately from the definition.

\begin{prop}\label{surj}
Let $\varphi : R_1 \rightarrow R_2$ be a homomorphism of rings. Then for any $x\in R_1$, the length of $x$ is greater than or equal to the length of $\varphi(x)$ in $R_2$. If $\varphi$ is an epimorphism, then $p(R_1)\geq p(R_2)$.
\end{prop}

\begin{thm}\cite[Proposition 4.5']{CLDR}\label{cldr}
Let $g(x,y) \in \R[x,y]$ be a polynomial of length $m$. Take a positive integer $r$ such that $2r>\deg_x g$. Then the polynomial $G(x,y)=g(x,y)(y-x^r)^2+1$ has length $m+1$.
\end{thm}
This theorem is a specialized version of \cite[Theorem 4.10]{CLDR} to the case $\R[x,y]$ (the much more general version can be applied to the various rings of type $A[x]$ for a commutative ring $A$, however we will not use it here). An immediate corollary is the following
\begin{cor}
The Pythagoras number of $\R[x,y]$ is infinite, i.e.
$$p(\R[x,y])=\infty.$$
\end{cor}

Note that the condition on $r$ above is not very restrictive, and we can inductively construct a sequence of polynomials $(F_m)_{m\geq 1}$ such that the length of $F_m$ is $m$ in the ring $\R[x,y]$.
\begin{df}\label{ciagi}
Let $f(x,y)$ be a strictly admissible polynomial. Consider a sequence of positive integers $(r_m)_{m\geq 1}$ such that
\begin{itemize}
\item $r_1>\deg f$ and $r_i> \sum_{k=1}^{i-1}r_{k}$
\item $f(x,x^{r_i})$ is a polynomial in $x$, which is either of odd degree, or it has even degree, and the leading coefficient is strictly positive. 
\end{itemize}

We define a sequence of polynomials in a following way:
\begin{itemize}
\item $F_1=1$
\item $F_n=F_{n-1}(y-x^{r_{n-1}})^2+1$

\end{itemize}
We call $(F_m)_{m\geq 1}$, a sequence of polynomials associated with $f(x,y)$.
\end{df}
The second condition follows from the admissibility of $f(x,y)$. Note that the condition $r_1>\deg f$ implies that the polynomial $f(x,x^{r_i})$ cannot be a constant polynomial for any $i$. By Theorem \ref{cldr}, length of each $F_n$ is equal to $n$ in the ring $\mathbb{R}[x,y]$.

Let now $f(x,y)$ be an admissible polynomial and consider the hypersurface $V$ given by the zero set of $z^2-f(x,y)$. Then the coordinate ring of $V$, $\R[V]=\R[x,y,z]/(z^2-f(x,y))$ is isomorphic to $\R[x,y,\sqrt{f(x,y)}]$, provided that $f(x,y)$ is not a square (if $f(x,y)$ is a square, then the former ring is isomorphic to $\R[x,y]\times \R[x,y]$ and the latter to $\R[x,y]$, both rings have infinite Pythagoras number). 

We have the following obvious lemma
\begin{lemma}\label{lemma}
Consider the polynomial $f(x,y)$ and $g(x,y)=(f\circ M)(x,y)$ for some invertible real matrix $M$. Then the rings $\R[x,y,z]/(z^2-f(x,y))$ and $\R[x,y,z]/(z^2-g(x,y))$ are isomorphic. \qed
\end{lemma}
At the end of this section, we provide sufficient condition for admissibility.
\begin{prop}Let $f(x,y)$ be a non-constant polynomial such that $f(0,0)=0$ and let $\deg f(x,y)=d$. Write $f(x,y)=g(x,y)+h(x,y)$, where $h(x,y)$ is the homogeneous part of degree $d$ and $\deg g(x,y)<d$. If $h(x,y)$ admits strictly positive value on $\R^2$, then $f(x,y)$ is an admissible polynomial.
\end{prop}
\noindent
\textit{Proof.}
Let $P\in \R^2$ be such that $h(P)>0$ and $M$ be the rotation matrix, such that it maps the halfline $OP$ to the positive part of $y$-axis. As $h$ is homogeneous, $h\circ M$ is strictly positive on the positive part of the $y$-axis, in other words, it contains a monomial $\alpha y^d$ with positive $\alpha$. Thus, $f\circ M$ is strictly admissible. \qed \\

The above condition is clearly not necessary, as $f(x,y)=-x^6+x^2y^2$ is admissible, but does not satisfy assumptions.

\begin{cor}
If $h(x,y)$ is a homogeneous admissible polynomial of degree $d$, then for any polynomial $g(x,y)$ with $\deg g(x,y) <d$, the polynomial $f(x,y):=g(x,y)-g(0,0)+h(x,y)$ is admissible.
\end{cor}

\section{Main result}

We will now focus on computing the Pythagoras number of the coordinate rings of a specific type of hypersurfaces. Our main result is as follows.

\begin{thm}\label{main}
Let $f(x,y)$ be an admissible polynomial which is not a square. Then 
$$p(\R[x,y, \sqrt{f(x,y)}])=\infty.$$
\end{thm}
\noindent
\textit{Proof.}
After a linear change of coordinates, $f(x,y)$ becomes a strictly admissible polynomial, hence, by the Lemma \ref{lemma}, we may assume that $f(x,y)$ is a strictly admissible polynomial. Let $(F_m)_{m \geq 1}$ be a sequence of polynomials associated with $f(x,y)$.
Assume to the contrary that the Pythagoras number of $\R[x,y, \sqrt{f(x,y)}]$ is finite and equal to a positive integer $L$, and consider the polynomial $F_{L+1}$. Every element $h \in \R[x,y, \sqrt{f(x,y)}]$ can be written in the form $h=h_1+\sqrt{f}h_2$, for some polynomials $h_1,h_2$.  By hypothesis, we have an equality
$$F_{L+1}=\sum_{i=1}^L (f_{i,1} + \sqrt{f}g_{i,1})^2$$
which translates into 
\begin{equation}\label{1}
F_{L+1}=\sum_{i=1}^L f^2_{i,1} + f\sum_{i=1}^Lg^2_{i,1}
\end{equation}
and 
$$\sum_{i=1}^Lf_{i,1}g_{i,1}\equiv 0,$$
for some polynomials $f_{i,1},g_{i,1} \in \R[x,y]$ (note that this last equality is not of interest to us).
After substitution $y=x^{r_{L}}$ we get
\begin{equation}\label{2}
F_{L+1}(x,x^{r_{L}})=\sum_{i=1}^L f^2_{i,1}(x,x^{r_{L}}) + f(x,x^{r_{L}})\sum_{i=1}^Lg^2_{i,1}(x,x^{r_{L}}).
\end{equation}
Recall that $F_{L+1}=F_L(y-x^{r_{L}})^2+1$, hence the above equation yields

\begin{equation}
1=\sum_{i=1}^L f^2_{i,1}(x,x^{r_{L}}) + f(x,x^{r_{L}})\sum_{i=1}^Lg^2_{i,1}(x,x^{r_{L}}).
\end{equation}
Since $f(x,y)$ is a strictly admissible polynomial, the constant term of the polynomial $f(x,x^{r_{L+1}})$ is zero. By construction of the associated sequence, $f(x,x^{r_{L}})$ is not a zero polynomial.  Moreover, the polynomial $ f(x,x^{r_{L}})\sum_{i=1}^Lg^2_{i,1}(x,x^{r_{L}})$ is either of odd degree, or of even degree with positive leading coefficient. Straightforward calculation shows that 
each $g_{i,1}$ is divisible by $y-x^{r_{L}}$ and each $f_{i,1} = a_i +(y-x^{r_{L}})f_{i,2}$ for $i=1,2,\dots, L$, where $a_i \in \mathbb{R}$ and $\sum_{i=1}^L a_i^2=1$. We may now apply orthogonal transformation over $\mathbb{R}[x,y]$ (see \cite{CLDR} or \cite[Theorem 8.1.2]{Prestel}), and assume that $f_{i,1}=(y-x^{r_{L}})f_{i,2}$ for $i=1,2,\dots, L-1$ and $f_{L,1}=1+(y-x^{r_{L}})f_{L,2}.$

We can rewrite equation (\ref{1}) as 
$$F_L(y-x^{r_{L}})^2+1=(y-x^{r_{L}})^2\sum_{i=1}^{L-1}f_{i,2}^2+1+2(y-x^{r_{L}})f_{L,2} +(y-x^{r_{L}})^2f_{L,2}^2 + (y-x^{r_{L}})^2f\sum_{i=1}^Lg^2_{i,2}.$$
After cancelling $1$'s we see that $(y-x^{r_{L}})$ divides $f_{L,2}$. If $f_{L,2}$
 is a nonzero polynomial, then the degree of the right hand side is at least $4r_{L}$. On the other hand, by construction, $2r_{L}$ is strictly larger than $ \deg F_L$, a contradiction. As a consequence, $f_{L,2}\equiv 0$ and we reduced the equation (\ref{1}) to
 $$F_{L}=\sum_{i=1}^{L-1} f^2_{i,2} + f\sum_{i=1}^Lg^2_{i,2}.$$
 With the above reasoning, we managed to reduce to number of polynomials $f_{i,1}$ by one. After repeating this procedure $L-2$ times, we obtain
 $$F_2=f_{1,L-1}^2 +f\sum_{i=1}^Lg^2_{i,L-1}$$
 which is equivalent to
 $$F_1(y-x^{r_1})^2+1=f_{1,L-1}^2 +f\sum_{i=1}^Lg^2_{i,L-1}.$$
After repeating the above procedure one last time, we get that $f_{1,L-1}\equiv 1 $ hence
$$F_1=f\sum_{i=1}^Lg^2_{i,L}.$$
This, however, is a contradiction as $F_1=1$ and $f$ is not a unit. This finishes the proof. \qed
 \newline

\begin{remark}
Denote by $\OO(V)$ the ring of regular functions \cite[Definition 3.2.1]{Bochnak} of the algebraic surface $V \subset \R^n$.
One can ask if the same reasoning would work for $\OO(V)$. However, the polynomials which allowed us to derive contradiction in the above proof were of the form $1+G(x,y)$, where $G(x,y)$ is a sum of squares. Every such element is a unit and a totally positive element of $\OO(V)$, hence by \cite[Theorem 7.3]{Mahe}, the length of any such polynomial is at most 4. Therefore, in order to compute $p(\OO(V))$, one has to use different methods. To the authors best knowledge $p(\OO(V))$ is not known for any surface.
\end{remark}
If we assume that $f(x,y)$ is a sum of squares we get the following:

\begin{thm}\label{sos}
Assume that $f(x,y)$ is a polynomial of length $n>1$. Then 
$$p(\R[x,y, \sqrt{f(x,y)}])=\infty.$$

\begin{proof}
Assume to the contrary, that $p(\R[x,y, \sqrt{f(x,y)}])=L$. Take $(G_m)_{m \in \N}$ to be any sequence of polynomials such that the length of $G_m$ is $m$ in $\R[x,y]$. In this case, $G_m$ can be written as 
$$G_m=\sum_{i=1}^L (f_i + \sqrt{f}g_i)^2$$
which translates into 
$$G_m=\sum_{i=1}^L f^2_i + f\sum_{i=1}^Lg^2_i$$
and $$\sum_{i=1}^Lf_{i}g_{i}\equiv 0,$$
for some polynomials $f_i,g_i \in \R[x,y]$.
We see, that the right hand side of the first equation contains at most $L+nL$ squares. Hence, by taking $m>L+nL$ we obtain a contradiction.
\end{proof}
\end{thm}
 We stress that we do not require $f(x,y)$ to be admissible.
As a final remark of this section. we would like to note that the above result cannot be generalized to the case where $f(x,y)$ is a strictly positive polynomial which is not a sum of squares. The existence of a zero of a polynomial is necessary for the proof of Theorem \ref{main}.

On the other hand, we have the following

\begin{thm}\label{-sos}
If $f(x,y)=- \sum_{i=1}^n f_i^2$, for some $f_i \in \R[x.y]$ then

$$
p(\R[x,y,z]/(z^2-f(x,y)))<\infty
$$
provided that $f(x,y)$ is a nonzero polynomial.
\end{thm}
\noindent
\textit{Proof.} Follows readily from \cite[Example 3.11]{CLDR}. \qed

\section{Applications and further problems}
In this last section we will provide some applications of Theorem \ref{main} and propose some further problems.

Consider the du Val singularities:
\begin{itemize}
\item $A_n: z^2 +x^2+y^{n+1}=0$, $n\geq 1$
\item $D_n: z^2+x^2y+y^{n-1}=0$, $n \geq 4$
\item $E_6: z^2+x^3+y^4=0$
\item $E_7: z^2 +x^3+xy^2=0$
\item $E_8: z^2+x^3+y^5=0$.
\end{itemize}
Note that the equations $A_{2k+1}$ describe a single point in $\R^3$, while the others are actual hypersurfaces in $\R^3.$ We may now state the following:
\begin{obs}
Let $Q$ be a one of equation from the above list. Then the affine $\R$-algebra $\R[x,y,z]/(Q)$ satisfies

$$p(\R[x,y,z]/(Q))
\begin{cases}
< +\infty  \ \ \ \text{for}  \ \ \ Q=A_{2k+1}\\
=+\infty \ \ \ \ \ \ \text{otherwise}.
\end{cases}
$$
\end{obs}

\noindent 
\textit{Proof.}
Any of the polynomials satisfying second part can be written in the form $z^2-f(x,y)$ where $f(x,y)$ is an admissible polynomial, hence the result follows from Theorem \ref{main}. The first part follows from Theorem \ref{-sos}. \qed \\

 In order to prove Theorem \ref{main} we had to assume that the polynomial $f(x,y)$ has a zero and is positive on a sufficiently large subset of $\R^2$. We propose the following

\begin{prob}
 Compute $$
p(\R[x,y,z]/(z^2-f(x,y))).
$$
in the following three cases:
\begin{itemize}
\item[a)] $f(x,y)$ is a strictly positive polynomial which is not a sum of squares,
\item[b)] $f(x,y)$ is a strictly negative polynomial such that $-f(x,y)$ is not a sum of squares,
\item[c)] $f(x,y)$ is indefinite, but it is not an admissible polynomial.
\end{itemize}
\end{prob}
 Consider the homogenized Motzkin polynomial $M(x,y,z)=z^6 +x^2y^4+x^4y^2-3x^2y^2z^2$. It is known that the cone of sums of squares of ternary sextics is strictly contained in the cone of positive semidefinite ternary sextics \cite{Reznick}. Since both cones are closed, there exists a positive $\epsilon>0$ real number such that the polynomial $M_1(x,y,z)=M(x,y,z)+\epsilon(x^6+y^6+z^6)$ is nonnegative, but it is not a sum of squares. Hence, after dehomogenizing, we see an interesting case in $a)$ and $b)$ above. For polynomials satisfying $c)$ see Example \ref{example} or any polynomial which is positive only on a bounded subset of $\R^2$.

We stress, that (up to translations) the above list together with polynomials satisfying Theorem \ref{main}, \ref{sos} or \ref{-sos} consists of the ring of all polynomials. We expect, that the Pythagoras number in $a)$ and $c)$ is infinite, and it is finite for polynomials satisfying $b)$. Based on the evidence shown in this paper, we formulate the following conjecture

\begin{conj}
Let $f(x,y,z) \in \R[x,y,z]$ be an irreducible polynomial. Then, the following conditions are equivalent:
\begin{itemize}
\item[a)]  $
p(\R[x,y,z]/(f(x,y,z))=+\infty$
\item[b)] the polynomial $f(x,y,z)$ is indefinite i.e. $f$ changes sign on $\R^3$ (equivalently, the ideal $(f(x,y,z))$ is real).

\end{itemize}
\end{conj}
Irreducibility condition is necessary in the above conjecture. For $f(x,y,z)=z^2$ the quotient ring has a surjection onto the ring $\R[x,y]$ hence by Proposition \ref{surj} its Pythagoras number is infinite.

The authors declare that no funds, grants, or other support were received during the preparation of this manuscript. \\

The authors have no relevant financial or non-financial interests to disclose.\\

Data sharing not applicable to this article as no datasets were generated or analysed during the current study.

\begin{small}

\vspace{5pt}
\noindent
Kacper Błachut

\noindent
Institute of Mathematics

\noindent
Faculty of Mathematics and Computer Science

\noindent
Jagiellonian University

\noindent
ul. Łojasiewicza 6, 30-348 Kraków, Poland

\noindent
e-mail: kacper.blachut@gmail.com

\vspace{5pt}

\noindent
Tomasz Kowalczyk

\noindent
Institute of Mathematics

\noindent
Faculty of Mathematics and Computer Science

\noindent
Jagiellonian University

\noindent
ul. Łojasiewicza 6, 30-348 Kraków, Poland

\noindent
e-mail: tomek.kowalczyk@uj.edu.pl

\end{small}

\end{document}